\documentclass[12pt]{amsart}
\usepackage[francais]{babel}
\usepackage{graphics}
\usepackage{graphicx}
\usepackage{pdfpages}

\title{Des crit\`eres de  transcendance inspir\'es par un texte de Kolberg dat\'e de 1962. II} 

\author{Labib Haddad}
\address{120 rue de Charonne, 75011 Paris, France}
\email{labib.haddad@wanadoo.fr}

\usepackage{amssymb}
\usepackage{amsmath}

\usepackage[T1]{fontenc}

\newcommand{\su}{\subsection*}
\newcommand{\head}{\section*}
\newcommand{\noi}{\noindent}

\newcommand{\Ž}{\'e}

\newcommand{\ˆ}{\`a}
\newcommand{\}{\`u}

\newcommand{\A}{\mathbb A}

\newcommand{\Q}{\mathbb Q}

\newcommand{\Z}{\mathbb Z}

\newcommand{\cal}{\mathcal}

\newcommand{\geqs}{\geqslant}

\newcommand{\ali} {\begin{aligned}}   
\newcommand{\ala} {\end{aligned}}

\newcommand{\stm}{\smallsetminus}

\newcommand{\inc}{\subset}

\newcommand{\bc}{\begin{cases}}
\newcommand{\ec}{\end{cases}}

\begin{document}
\maketitle

\thispagestyle{empty}

\markboth{Labib Haddad}{Oddmund Kolberg}

Dans un article de 1962, voir [1], {\sc  O. Kolberg} \Žnonce et \Žtablit ceci.

\su{Th\Žor\me} On se donne un entier $a>0$, un nombre rationnel $r$, distinct de $-1,-2,\dots, -a+1$, et un polyn\™me $P(z)$ non nul, \ˆ coefficients alg\Žbriques.  Alors, pour  $x$ alg\Žbrique tel que $0 < |x| < 1/e$, la somme de la s\Žrie suivante, $S$, est un nombre transcendant :
$$S = \sum^\infty_{n=1}\frac{(n+r)^{n-a}P(n)}{n!} x^n.$$
Ce r\Žsultat est-il aussi isol\Ž qu'il le para\"t, s'appliquant \ˆ une  famille assez particuli\re de s\Žries enti\res ?

\noi Dans la premi\re partie de ce texte, voir [0] en bibliographie et errata,  on avait donn\Ž un exemple du m\me ordre pouvant tenir compagnie \ˆ celui de Kolberg.

\su{Exemple} Pour tout entier $a\in\Z$, tout polyn\™me $P(z)$ non  nul \ˆ coefficients alg\Žbriques, tout nombre  alg\Žbrique $x$ tel que $0 < |x| < 1/e$, la somme de la s\Žrie
$$\sum_{n=2}^\infty (n-1)n^{n+a}  P\left(\frac{1}{n}\right) \frac{x^n}{n!}$$
est un nombre transcendant ! 

\

\noi Cependant, \ˆ y regarder de pr\s, ces s\Žries sont de la m\me famille que celles de Kolberg !

\noi En effet, si le degr\Ž du polyn\™me $P(z)$ est \Žgal \ˆ $d$, alors 
$$(n-1)n^dP\left(\frac{1}{n}\right)  = Q(n)$$
est un polyn\™me en $n$, et la s\Žrie s'\Žcrit
$$\sum_{n=2}^\infty n^{n+a-d}  Q(n) \frac{x^n}{n!}$$

\

\noi Cet exemple mentionn\Ž  dans la premi\re partie de notre texte, voir [0],  \Žtait donn\Ž comme illustration  de {\it la kolbergisation} ! Lequel proc\Žd\Ž est  une {\bf imitation} de celui qu'utilise Kolberg pour la d\Žmonstration de son th\Žor\me !

\

\noi On se reportera \ˆ cette premi\re partie du texte, [0], pour tous les d\Žtails au sujet de la kolbergisation. On va en rappeler ici, tout de m\me, succinctement, l'essentiel.

\su{Kolbergisation} Ce proc\Žd\Ž utilise les ingr\Ždients suivants. Un changement de variable, $x = te^{-t}$. Le th\Žor\me de Lindemann.  La notion de s\Žrie associ\Že. Un crit\re de transcendance pour des fractions rationnelles, et enfin les quatuors. 

\noi Voici quelques rappels.

\

\head{Le th\Žor\me de Lindemann}

Pour $t$ alg\Žbrique non nul, $x = te^{-t}$ est transcendant. Donc, lorsque $x$ est alg\Žbrique non nul, $t$ est transcendant.

\

\head{Le crit\re de transcendance}

\ 

Soit $g(s)$ une fraction rationnelle en $s$  dont les coefficients sont des nombres alg\Žbriques. On d\Žsigne par $E$ l'ensemble (exceptionnel) des nombres entiers $n\in \Z$ pour lesquels la fraction rationnelle $s^ng(s)$ est constante. On se donne un nombre rationnel $r \in \Q\stm E$. Alors, si  $t$ est un nombre transcendant, le nombre $t^rg(r)$ l'est  \Žgalement.

\

\noi Voir la d\Žmonstration dans [0].

\

\head{La notion de s\Žrie associ\Že}

\
 
\noi On se donne une s\Žrie enti\re en $x$
$$H(x) = \sum_{n=0}^\infty u_n\frac{x^n}{n!}.$$
On d\Žfinit la {\bf s\Žrie associ\Že}
$$G(t) = H(te^{-t}) = \sum_{n=0}^\infty u_n\frac{t^ne^{-nt}}{n!}= \sum_{n=0}^\infty v_k\frac{t^k}{k!}.$$
On \Žtablit alors ceci. On a $u_0 = v_0$ et, pour  $n\geqs 1$,
$$\boxed{u_n = v_n +  \dots + \binom {n-1}{m-1} n^{n-m} v_m + \dots+ n^{n-1}v_1,}$$
ce qui \Žquivaut encore \ˆ :
$$\boxed{v_n = u_n  + \dots+  \binom n m (-m)^{n-m} u_m + \dots  +    (-1)^{n-1}nu_1.}$$

\

\noi Voir les d\Žmonstrations dans [0]. [{\color {red} Noter que, dans [0], il est \Žcrit, par erreur, $(-n)^{n-m}$ au lieu de $(-m)^{n-m}$}.]

\

\

\head{Les quatuors}

\

\noi {\bf Une suite} est, par convention, une famille ind\Žx\Že par l'ensemble $\Z$ des entiers.

\

\noi {\bf Un quatuor} est form\Ž de quatre suites de fonctions, $(F,G,H,K)$ :
$$K=(K_k(x,y))_{k\in \Z}, H  = (H_k(x,y))_{k\in \Z},$$
$$F = (F_k(t,y))_{k\in\Z}, G= (G_k(t,y))_{k\in\Z},$$
li\Žes par les relations suivantes
$$K_k(x,y) = x^yH_k(x,y)$$
$$G_k(t,y) = H_k(te^{-t},y)$$
$$F_k(t,y) = K_k(te^{-t},y)$$ 
et satisfaisant les condition suivantes :
\[x\frac{d}{dx}K_{k+1}(x,y) = K_k(x,y)\tag*{(1k)}\]
Il s'ensuit quelques autres conditions du m\me genre, portant sur les fonctions $H_k, G_k, F_k,$ \Žgalement. Voici une liste de toutes ces relations au sein des quatuors, d\Žtaill\Žes dans [0].

\

\

  $\bc

K_k(x,y) = x^yH_k(x,y)\\ 
\\

G_k(t,y) = H_k(te^{-t}, y)\\ 
\\

F_k(t,y) = K_k(te^{-t},y) =  t^ye^{-yt}G_k(t,y)\\
\\

\displaystyle
(1k) \ x\frac{d}{dx}K_{k+1}(x,y) = K_k(x,y)\\
\\

\displaystyle
(2k) \ K_{k+1}(x,y) = \int_0^x \frac{K_k(z,y)}{z}dz\\

\\
\displaystyle
(3k) \ x\frac{d}{dx} H_{k+1}(x,y) + yH_{k+1}(x,y) = H_k(x,y)\\
\\

\displaystyle
(4k) \ x^yH_{k+1}(x,y) =  \int_0^x  z^{y-1}H_k(z,y)dz \\
\\

\displaystyle
(5k) \ \frac{t}{1-t}\frac{d}{dt} G_{k+1}(t,y) + yG_{k+1}(t,y) = G_k(t,y)\\
\\

\displaystyle
(6k) \ t^ye^{-yt}G_{k+1}(t,y) = \int_0^t \frac{1-z}{z}z^ye^{-yz}G_k(z,y)dz\\
\\

\displaystyle
(7k) \ \frac{d}{dt}F_{k+1}(t,y) = \frac{1-t}{t}F_k(t,y)\\
\\

\displaystyle
(8k) \ F_{k+1}(t,y) = \int_0^t \frac{1-z}{z} F_k(z,y) dz\\
\\
 
\ec$ 

\

\

\su{Op\Žrations sur les quatuors} De mani\re naturelle, on d\Žfinit  le d\Žcal\Ž d'odre $d$ du quatuor  $\cal Q$ comme \Žtant le quatuor $\cal Q^{\to d}=$
$$K^{\to d}=(K_{k+d}(x,y))_{k\in \Z}, H^{\to d}  = (H_{k+d}(x,y))_{k\in \Z},$$
$$F^{\to d} = (F_{k+d}(t,y))_{k\in\Z}, G^{\to d}= (G_{k+d}(t,y))_{k\in\Z}.$$
De m\me, on d\Žfinit de mani\re naturelle les combinaisons lin\Žaires de quatuors, $\lambda_1\cal Q^1 +\dots \lambda_p\cal Q^p$, qui sont elles-m\mes des quatuors. Voir d\Žtails dans la premi\re partie, [0].
\

\

\su{Le quatuor de Kolberg} C'est le quatuor $\cal K = (F,G,H,K)$  o\ 
$$H_k(x,y) = T_k(x,y) = \sum^\infty_{n=0} (y+n)^{n-k}\frac{x^n}{n!}.$$ 
$$K_k(x,y) = x^yT_k(x,y) = f_k(x,y).$$

\

\head{Description de la kolbergisation}

\

\centerline {\bf Voici en bref ce qu'est la kolbergisation !}

\

On se donne un quatuor, {\it fertile}, $\cal K = (F,G,H,K)$, c'est-\ˆ-dire tel que chacune des fonctions $F_k(t,y)$ soit de la forme
\[F_k(t,y) = t^y R_k(t,y), \ R_k(t,y) \in  \A(y)(t)^*.\tag*{(ad hoc)}\]
o\ $\A$ d\Žsigne le corps des nombres alg\Žbriques. Autrement dit, $R_k(t,y)$ est une fraction rationnelle (non nulle) en $t$ dont les coefficients sont des fractions rationnelles en $y$ \ˆ coefficients alg\Žbriques ! 

\

\noi Pour toute partie $I\inc \Z$, on d\Žsigne par $Y(I)$ l'ensemble de tous les p\™les  de toutes ces  fractions rationnelles en $y$ qui interviennent comme coefficients dans les $F_k(t,y), k\in I$.

\

\noi {\it Soient $A_k, k\in I$, des nombres alg\Žbriques non tous nuls. On se donne un nombbre $r \in \Q\stm Y(I)$. La combinaison lin\Žaire 
$$L = \sum_{k\in I} A_kK_k(x,r),$$
est un nombre transcendant lorsque $x$ est alg\Žbrique non nul} \ !

\su{Le raisonnement est le suivant} On utilise le changement de variable, $x = te^{-t}$. Il vient $L = \sum_k A_kF_k(t,y)= t^rg(t)$  auquel le crit\re de transcendance s'applique. Ainsi, pour $x$ alg\Žbrique, on sait que $t$ est transcendant (Lindemann) donc $L$ est transcendant. \qed

\

\noi Le {\bf quatuor fertile} le plus g\Žn\Žral est engendr\Ž par une fonction quelconque,  $F_0(t,y)$ non nulle, de la forme 
$$F_0(t,y) = t^y R_0(t,y), \ R_0(t,y)\in \A(y)(t).$$
Les autres fonctions $F_k(t,y)$ s'obtiennent alors, en avant et en arri\re, \ˆ l'aide des formules (7k) et (8k) :
$$F_k(t,y) = \frac{t}{1-t}\frac{d}{dt}F_{k+1}(t,y) , \ F_{k+1}(t,y) = \int_0^t \frac{1-z}{z} F_k(z,y) dz.$$
Si $F_{k+1}(t,y)$ a la forme ad hoc, la fonction $F_k(t,y)$ l'aura aussi, par d\Žrivation. Il faut encore que l'int\Žgration redonne la forme ad hoc. C'est la seule chose \ˆ v\Žrifier pour savoir si le quatuor est vraiment fertile !

\

\centerline{\bf Voil\ˆ r\Žsum\Že la kolbergisation !}

\

\su{Remarque} Tout d\Žcal\Ž d'un quatuor fertile est, bien entendu, fertile. D'autre part, toute combinaison lin\Žaire {\bf non nulle}, \ˆ coefficients alg\Žriques, de quatuors fertiles est un quatuor fertile !

\head{Illustration \\ \small Un exemple de kolbergisation}

\

On d\Žsigne par $\cal K^\sharp= (F,G,H,K)$ le quatuor engendr\Ž par
$$G_0(t,y) = \left(1+\frac{2}{y}+t^2\right)e^{yt},$$
autrement dit,
$$F_0(t,y) = \left(1+\frac{2}{y}+t^2\right)t^y.$$
En it\Žrant  l'int\Žgration
$$F_{k+1}(t,y) = \int_0^t \frac{1-z}{z} F_k(z,y) dz,$$
on v\Žrifie que le quatuor est fertile !

\

\noi D'autre part, on 
$$G_0(t,y) = \left(1+\frac{2}{y}+t^2\right)e^{yt}= \frac{y+2}{y} + \sum_{n=1}^\infty (y^2 + 2y + n(n-1))y^{n-2}\frac{t^n}{n!}.$$
La s\Žrie $G_0(t,y)$ est l'associ\Že de la s\Žrie  :
$$H_0(x,y) = \frac{y+2}{y} + \sum_{n=1}^\infty u_n\frac{x^n}{n!}= \sum_{n=0}^\infty u_n\frac{x^n}{n!}.$$
On a $u_0 = \displaystyle\frac{y+2}{y}$ et,
pour $n\geqs 1$, les $u_n$ s'obtiennent par les formules suivantes :
$$u_n = \sum_{m=1}^n \binom {n-1} {m-1} n^{n-m}(y^2 + 2y + m(m-1))y^{m-2}.$$ 
On v\Žrifie alors que, pour tout $n\geqs 0$, on a
$$u_n = (y+2)(y^2+ 2ny + 2n^2- n)(y+n)^{n-3}.$$

\

$u_0 = (y+2)/y$ 

$u_1 =y+2$

$u_2 = y^2 + 4y +6$

$u_3 = (y+2)(y^2 + \ 6y \ + 15)$
                            
$u_4  = (y+2)(y^2+ \ 8y \ + 28)(y+ 4)$

$u_5  = (y+2)(y^2 + 10y + 45)(y+5)^2$

$u_6 = (y+2)(y^2 + 12y + 66)(y+6)^3$

 $u_7 =  (y+2)(y^2+ 14 y + 91)(y+7)^4$.

$$H_0(x,y) =  \frac{y+2}{y} + \sum_{n=1}^\infty (y+2)(y^2+ 2ny + 2n^2 n)(y+n)^{n-3}\frac{x^n}{n!},$$
$$K_0(x,y) =  x^yH_0(x,y) = \sum_{n=0}^\infty (y+2)(y^2+ 2ny + 2n^2 n)(y+n)^{n-3}\frac{x^{n+y}}{n!}.$$
\`A l'aide des formules (1k) et (2k), qui fournissent $K_k(x,y)$ en fonction de $K_{k+1}(x,y)$ et vice versa, on obtient
$$K_k(x,y) =  x^y\sum_{n=0}^\infty (y+2)(y^2+ 2ny + 2n^2- n)(y+n)^{n-k-3}\frac{x^n}{n!}.$$
Aux fins de comparaison, les fonctions correspondantes dans le quatuor de Kolberg sont les
$$f_k(x,y) = x^y\sum^\infty_{n=0} (y+n)^{n-k}\frac{x^n}{n!}.$$
On peut alors \Žnoncer le r\Žsultat suivant :

\

\noi {\it On se donnne un nombre rationnnel $r\neq 0, -1, -2,\dots$, et  des  nombres alg\Žbriques  $A_k$. Toute combinaison lin\Žaire non nulle,
$$L(x) = \sum_k A_k x^r\sum_{n=0}^\infty (r+2)(r^2+ 2nr + 2n^2 -n)(r+n)^{n-k-3}\frac{x^n}{n!},$$
est un nombre transcendant, pour tout $x$ alg\Žbrique, $0 < |x| < 1/e$}.

\

On en d\Žduit ceci qui tient compagnie au th\Žor\me  de Kolberg.

\

\su{Th\Žor\me di\se} On se donne un entier $a>0$, un polyn\™me non nul, $P(z)$, \ˆ coefficients   alg\Žbriques, un nombre rationnel $r \neq -1, -2, \dots$.  Alors, pour  $x$ alg\Žbrique tel que $0 < |x| < 1/e$, la somme $S$ de la s\Žrie suivante est un nombre transcendant :
$$S = \sum^\infty_{n=0}\frac{(r^2+ 2nr + 2n^2 -n)(n+r)^{n-a}P(n)}{n!} x^n.$$

\

\

\noi {\bf Pour  conclure}, on dira que la kolbergisation peut fournir, sans doute, encore d'autres r\Žsultats semblables pour d'autres classes de s\Žries enti\res \ˆ valeurs transcendantes. Cela m\Žrite d'\tre \Žlabor\Ž !

\

\

\head{Bibliographie}

\

\noi [0] {\sc L. Haddad}, {\sl Des crit\`eres de  transcendance inspir\'es par un texte de Kolberg dat\'e de 1962}, {\tt arXiv\!:2103.13307v1 [math.NT] 22 Mar 2021.}

\

\noi ERRATUM :  dans la formule de la page 6 de cette note, [0], qui donne  les $v_n$ en fonction des $u_m$, il faut remplacer $(-n)^{n-m}$ par  $(-m)^{n-m}$.

\

\noi [1] {\sc O. Kolberg}, {\sl A class of power series with transcendental sums for algebraic values of the variable}. \AA rbok Univ. Bergen, Mat.-Naturv. Ser. (1962) No. 18, 6 p.

\

\

\

\enddocument